\documentstyle{amsppt}
\magnification1200
\pagewidth{6.5 true in}
\pageheight{9.25 true in}
\NoBlackBoxes

\def\phi{\varphi}

\def\cbar{\overline{\chi}}

\topmatter
\title Extreme values of $|\zeta(1+it)|$
\endtitle
\author Andrew Granville and K. Soundararajan
\endauthor
\address{D{\'e}partment  de Math{\'e}matiques et Statistique,
Universit{\'e} de Montr{\'e}al, CP 6128 succ Centre-Ville,
Montr{\'e}al, QC  H3C 3J7, Canada}\endaddress
\email{andrew{\@}dms.umontreal.ca}
\endemail
\address{Department of Mathematics, University of Michigan, Ann Arbor,
Michigan 48109, USA} \endaddress
\email{ksound{\@}umich.edu} \endemail
\thanks{Le premier auteur est partiellement soutenu par une bourse
de la Conseil  de recherches en
sciences naturelles et en g\' enie du Canada. The second  author is
partially supported by the National Science Foundation.}
\endthanks
\dedicatory 
To Professor K. Ramachandra on the occasion of his seventieth 
birthday 
\enddedicatory
\endtopmatter

\document

\head 1. Introduction \endhead

\noindent
 Improving on a result of J.E. Littlewood, N. Levinson [3] showed that there
are arbitrarily large $t$ for which $|\zeta(1+it)| \ge e^{\gamma}
\log_2 t +O(1)$. (Throughout $\zeta(s)$ is the Riemann-zeta function,
and $\log_j$ denotes the $j$-th iterated logarithm, so that $\log_1n
=\log n$ and $\log_j n=\log (\log_{j-1}n) $ for each $j\geq 2$.)
The best upper bound known is Vinogradov's
$|\zeta(1+it)| \ll (\log t)^{2/3}$.

Littlewood had shown that  $|\zeta(1+it)|\lesssim 2 e^{\gamma} \log_2 t$
assuming the Riemann Hypothesis, in fact by showing that
the value of $|\zeta(1+it)|$ could be closely approximated by
its Euler product for primes up to $\log^2 (2+|t|)$ under this assumption.
Under the further hypothesis
that  the Euler product up to $\log (2+|t|)$ 
still serves as a good approximation,
Littlewood conjectured that $\max_{|t|\leq T} |\zeta(1+it)| \sim e^{\gamma}
\log_2 T$, though later he wrote in [5] (in 
connection with a $q$-analogue): ``{\sl there is perhaps no good
reason for believing ... this hypothesis}''.

Our
Theorem 1 evaluates the frequency with which such extreme values are
attained; and if this density function were to persist to the end of
the viable range then this implies the conjecture that
$$
\max_{t\in [T,2T]} |\zeta(1+it)| = e^{\gamma} (\log_2 T+\log_3 T
+C_1+o(1)), \tag{1.1a}
$$
for some constant $C_1$. In fact it may be that $C_1=C+1-\log 2$, where
$$
C = \int_0^2 \log I_0(t) \frac{dt}{t^2} + \int_2^\infty (\log
I_0(t)-t) \frac{dt}{t^2} ,
$$
and $I_0(t):={\Bbb E}(e^{\text{Re} (tX)}) = \sum_{n=0}^{\infty}
(t/2)^{2n}/n!^2$ is the Bessel function (with $X$  a random variable
equidistributed on the unit circle).  In Theorem 2 we show that there are
arbitrarily large $t$ for which $|\zeta(1+it)| \ge e^{\gamma}(\log_2
t + \log_3 t- \log_4 t +O(1))$, which improves upon Levinson's
result but falls a little short of our conjecture.

Levinson  also showed that $1/|\zeta(1+it)| \ge
\frac{6e^{\gamma}}{\pi^2} (\log_2 t -\log_3 t+ O(1))$ for
arbitrarily large $t$.  Theorem 1 exhibits even smaller values of
$|\zeta(1+it)|$ and determines their frequency.  Extrapolating
Theorem 1 we are also led to conjecture that
$$
\max_{t\in [T,2T]} 1/|\zeta(1+it)| = \frac{6e^{\gamma}}{\pi^2}
     (\log_2 T+\log_3 T +C_1+o(1)); \tag{1.1b}
$$
but only succeed in proving that $1/|\zeta(1+it)| \ge
\frac{6e^{\gamma}}{\pi^2} (\log_2 t - O(1))$ for arbitrarily large
$t$.  K. Ramachandra [6] has obtained results 
analogus to Levinson's in short intervals, 
and R. Balasubramanian, Ramachandra and A. Sankaranarayanan 
[1] have considered extreme values of $|\zeta(1+it)|^{e^{i\theta}}$ 
for any $\theta \in [0,2\pi)$.

To be more precise let us define, for $T, \tau\geq 1$,
$$
\align  {\Phi}_T(\tau):&= \frac{1}{T} \text{meas} \{ t\in [T,2T]: \
\ |\zeta(1+it)| >e^{\gamma}\tau\}, \\
\text{and} \ \  {\Psi}_T(\tau):&= \frac{1}{T} \text{meas} \{ t\in
[T,2T]: \ \ |\zeta(1+it)| <\tfrac{\pi^2}{6e^{\gamma}\tau}\}.
\endalign
$$

\proclaim{Theorem 1}  Let $T$ be large.  Uniformly in the range
$1\ll \tau \le \log_2 T -20$ we have
$$
{\Phi}_T(\tau) = \exp\Big(-\frac{2e^{\tau-C-1}}{\tau} \Big(1 +
O\Big(\frac{1}{\tau^{\frac 12}} +\Big(\frac{e^{\tau}}{\log T}\Big)^{\frac 12}
 \Big)\Big)\Big) ,
$$
where $c$ is a positive constant. The same asymptotic also holds for
$\Psi_T(\tau)$.
\endproclaim

With a judicious application of the pigeonhole principle we can
exhibit even larger values of $|\zeta(1+it)|$, indeed of almost the
same quality as the conjectured (1.1a).

\proclaim{Theorem 2}  For large $T$ the subset of points $t \in
[0,T]$ such that
$$
|\zeta(1+it)| \ge e^{\gamma}(\log_2 T + \log_3 T- \log_4 T- \log
A+O(1))
$$
has measure at least $T^{1-\frac 1A}$, uniformly for any $A\geq 10$.
\endproclaim

One can also establish results analogous to Theorems 1 and 2 for the
distribution of values of $|L(1,\chi)|$ where $\chi$ ranges over all
non-trivial characters modulo a large prime $p$ (see section 7 for
further details). In fact Theorems 1 and 2 hold almost verbatim,
just changing $T$ to $p$. If one also averages over $p$ in a dyadic
interval $P\le p\le 2P$ then one can obtain asymptotics for the
distribution function in the wider range $1\ll \tau \leq \log_2 P +\log_3 P -
O(1)$ (which we expect is the full range, up to the explicit value of
the ``$O(1)$'').

As in [2] we can compare the distribution of $\zeta(1+it)$ with that of an
appropriate probabilistic model.  Let $X(p)$ denote independent 
random variables uniformly distributed on the unit circle, for each 
prime $p$.
We extend $X$ multiplicatively to all integers $n$: that is 
set $X(n)=\prod_{p^{\alpha} \parallel n} X(p)^{\alpha}$. 
We wish to compare the distribution of values of $\zeta(1+it)$ with 
the distribution of values of the random Euler products $L(1,X):= 
\prod_p (1-X(p)/p)^{-1}$ (these products converge with probability $1$).  
Now define
$$
\Phi(\tau) = \text{Prob}(|L(1,X)| \ge e^{\gamma} \tau) \ \ \text{and} \ \
\Psi(\tau) = \text{Prob} (|L(1,X)| \le \tfrac{\pi^2}{6 e^{\gamma }\tau} ). 
$$
By the same methods one can show that $\Phi(\tau)$ and $\Psi(\tau)$ satisfy
the same asymptotic as $\Phi_T(\tau)$ as in Theorem 1, but for arbitrary $\tau$
(see the remarks immediately after the
proof of Theorem 1).

\head 2.  Preliminaries \endhead

\noindent We collect here some standard facts on $\zeta(s)$ which
will be used later.

\proclaim{Lemma 2.1}  Let $y\ge 2$ and $|t| \ge y+3$ be real numbers.
Let $\frac 12 \le \sigma_0 <1$ and suppose that the
rectangle $\{ z: \ \ \sigma_0 <\text{Re}(z) \le 1,
\ \ |\text{Im}(z) -t| \le y+2\}$ is free of zeros of $\zeta(z)$.
Then for any $\sigma_0 < \sigma \le 2$ and $|\xi-t|\le y$ we
have
$$
|\log \zeta(\sigma +i\xi)| \ll \log |t|  \log (e/(\sigma -\sigma_0)).
$$
Further for $\sigma_0 < \sigma \le 1$ we have
$$
\log \zeta(\sigma+it)= \sum_{n=2}^{y}
\frac{\Lambda(n)}{n^{\sigma+it} \log n} + O\Big(
\frac{\log |t|}{(\sigma_1-\sigma_0)^2}y^{\sigma_1-\sigma}\Big),
$$
where we put $\sigma_1 = \min(\sigma_0+\frac{1}{\log y},
\frac{\sigma+\sigma_0}{2})$.
\endproclaim
\demo{Proof}  The first assertion follows from Theorem 9.6(B) of
Titchmarsh [8].  In proving the second assertion we may plainly
suppose that $y \in {\Bbb Z}+\frac 12$.  Then Perron's formula
gives, with $c= 1-\sigma+\frac{1}{\log y}$,
$$
\align
\frac{1}{2\pi i} \int_{c-iy}^{c+iy} \log \zeta(\sigma+it+w)\frac{y^w}{w}
dw &= \sum_{n=2}^{y}
\frac{\Lambda(n)}{n^{\sigma+it} \log n} + O\Big(\frac{1}{y}
\sum_{n=1}^{\infty} \frac{y^c}{n^{\sigma+c}} \frac{1}{|\log (y/n)|}\Big) \\
&= \sum_{n=2}^{y}
\frac{\Lambda(n)}{n^{\sigma+it} \log n} + O(y^{-\sigma} \log y).
\tag{2.1}\\
\endalign
$$
We now move the line of integration to the line Re$(w)=\sigma_1
-\sigma <0$.  Our hypothesis ensures that the integrand is regular
over the region where the line is moved, except for a simple pole at
$w=0$ which leaves the residue $\log \zeta(\sigma+it)$. Thus the
left side of (2.1) equals $\log \zeta(\sigma+it)$ plus
$$
\frac{1}{2\pi i} \Big( \int_{c-iy}^{\sigma_1-\sigma-iy}
+\int_{\sigma_1-\sigma -iy}^{\sigma_1-\sigma+iy}
+\int_{\sigma_1-\sigma +iy}^{c+iy}\Big) \log \zeta(\sigma+it+w)\frac{y^w}{w}
dw
\ll \frac{\log |t|}{(\sigma_1-\sigma_0)^2}y^{\sigma_1-\sigma},
$$
upon using the first part of the Lemma.

\enddemo

Using Lemma 2.1 we shall show that most of the
time we may approximate $\zeta(s)$ by a short Euler product.

\proclaim{Lemma 2.2}  Let $\frac 12 < \sigma \le 1$ be
fixed and let $T$ be large.  Let $T/2 \ge y \ge 3$ be a real number.
The asymptotic
$$
\log \zeta(\sigma+it) = \sum_{n=2}^{y}
\frac{\Lambda(n)}{n^{\sigma+it} \log n}
+ O( y^{(\frac 12-\sigma)/2} \log^3 T)
$$
holds for all $t\in (T,2T)$ except for a set of measure
$\ll T^{5/4-\sigma/2} y (\log T)^5$.
\endproclaim

\demo{Proof}  This follows upon using the zero-density result
$N(\sigma_0,T) \ll T^{3/2-\sigma_0} (\log T)^5$ (see Theorem 9.19 A
of [8]) and appealing to Lemma 2.1 (taking $\sigma_0 = (1/2
+\sigma)/2$ there).
\enddemo

\head 3.  Approximating $\zeta(1+it)$ by a short Euler product \endhead

\proclaim{Lemma 3.1}  Suppose $2\le y\le z$ are real numbers.
Then for arbitrary complex numbers $x(p)$ we have
$$
\frac{1}{T} \int_{T}^{2T} \Big| \sum_{y\le p\le z} \frac{x(p)}{p^{it}}
\Big|^{2k} \ll \Big(k\sum_{y\le p \le z} |x(p)|^2 \Big)^k
+ T^{-\frac 23} \Big( \sum_{y\le p\le z} |x(p)|\Big)^{2k}
$$
for all integers $1\le k \le \log T/(3\log z)$.
\endproclaim
\demo{Proof}  The quantity we seek to estimate is
$$
\sum\Sb p_1, \ldots , p_k \\ y\le p_j \le z\endSb
\sum\Sb q_1, \ldots, q_k \\ y\le q_j \le z\endSb
\overline{x(p_1) \cdots x(p_k)} x(q_1) \cdots x(q_k)
\frac{1}{T} \int_{T}^{2T}
\Big(\frac{p_1\cdots p_k}{q_1 \cdots q_k}\Big)^{it} dt.
$$
The diagonal terms $p_1\cdots p_k = q_1 \cdots q_k$ contribute
$$
\ll k! \Big( \sum_{y\le p\le z} |x(p)|^2 \Big)^k.
$$
If $p_1\cdots p_k \neq q_1 \cdots q_k$ then as both quantities are
below $z^k \le T^{\frac 13}$ we have that
$$
\frac{1}{T}\int_{T}^{2T}
\Big(\frac{p_1\cdots p_k}{q_1 \cdots q_k}\Big)^{it} dt
\ll \frac{1}{T |\log (p_1\cdots p_k/q_1\cdots q_k)|}
\ll T^{-\frac 23}.
$$
Hence the off diagonal terms contribute
$\ll T^{-\frac 23}(\sum_{y\le p\le z} |x(p)|)^{2k}$,
proving the Lemma.

\enddemo

Define $\zeta(s;y):=\prod_{p\leq y}  (1-p^{-s})^{-1}$.

\proclaim{Proposition 3.2}  Let $T$ be large and let $\log T (\log_2
T)^4 \ge y \ge e^2 \log T$ be a real number. Then there is a
positive constant $c$ such that
$$
\zeta(1+it) = \zeta(1+it;y) \Big( 1+ O\Big( \frac{\sqrt{\log T}}{\sqrt{y}
 \log_2 T} \Big)\Big)
$$
for all $t\in (T,2T)$ except for a set of measure at most
$T \exp(-\log T/50\log_2 T)$.
\endproclaim
\demo{Proof} Setting $z= (\log T)^{100}$ we deduce from Lemma 2.2
that $\zeta(1+it) = \zeta(1+it;z) (1+O(1/\log T))$ for
all $t\in (T,2T)$ except for a set of measure at most $T^{4/5}$.
Using Lemma 3.1 with $k=[\log T/(300 \log_2 T)]$ and $x(p)=1/p$
we get that
$$
\align \frac{1}{T} \int_{T}^{2T} \Big|\sum_{y\le p\le z}
\frac{1}{p^{1+it}} \Big|^{2k} dt &\ll \Big(k \sum_{y\le p\le z}
\frac{1}{p^2}\Big)^k + T^{-\frac 23}
\Big(\sum_{y\le p\le z} \frac 1p\Big)^{2k} \\
&\ll \Big( \frac{\log T}y\Big)^k \Big(\frac{1}{10 \log y}\Big)^{2k}
+ T^{-\frac 13},\\
\endalign
$$
and so 
$$
\Big|\sum_{y\le p\le z} \frac{1}{p^{1+it}} \Big| \le \frac{\sqrt{\log T}}{
\sqrt{y} \log y}
$$
for all $t\in [T,2T]$ except for a set of measure $\le T\exp(-\log T/49\log_2 T)$.
The Proposition thus follows, by combining the above estimates, since 
$$
\zeta(1+it;y)=\zeta(1+it;z) \exp\Big( -\sum_{y\le p\le z}
\Big(\frac{1}{p^{1+it}} + O\Big(\frac{1}{p^2}\Big)\Big)\Big) .
$$
\enddemo

\head 4.  Moments of short Euler products \endhead

\noindent In this section we show how to evaluate large
moments of the short Euler products obtained in \S 3.

\proclaim{Theorem 4.1} Let $ \log T(\log_2 T)^4 \ge y \ge e^2 \log T$
be a real number. Let $z=\delta k$ where $\delta = \pm 1$ and
$2\le k \le \log T/(e^{10}\log (y/\log T))$ is an integer.
Then
$$
\align
\frac 1T \int_T^{2T} |\zeta(1+it;y)|^{2z} &=
\sum\Sb n=1\\ p|n \implies p\le y\endSb^{\infty}
\frac{d_z(n)^2}{n^2} \Big( 1+ O\Big(\exp\Big(-\frac{\log T}{2(\log_2 T)^4}
\Big)\Big)\Big)\\
&=\prod_{p\le k} \Big(1-\frac{\delta}{p}\Big)^{-2k\delta}
\exp\Big(\frac{2k}{\log k}\Big(C+ O\Big(\frac{k}{y} +\frac{1}{\log k}\Big)
\Big)\Big).\\
\endalign
$$
\endproclaim

Throughout this section let $z$, $y$, $k$, $\delta$ be
as in Theorem 4.1.  If $k\le 10^6$ then we divide $[1,y]$ into
the intervals $I_0= [k,y]$ and $I_1=[1,k)$ and take here $J:=1$.
If $k>10^6$ then we define $J:=[4\log_2 k/\log 2]+1$ and
divide $[1,y]$ into the $J+1$-intervals $I_0=[k,y]$,
$I_j = [k/2^j,k/2^{j-1})$ for $1\le j\le J-1$,
and $I_J=[1,k/2^J) \subset [1,k/(\log k)^4]$.
Given a subset $R$ of the index set $\{0,\ 1,\ \ldots, \ J\}$ we
define ${\Cal S}(R)$ to be the set of integers $n$ whose prime
factors all lie in $\cup_{r\in R} I_r$.  We also define
$$
\zeta(s;R):= \prod_{p\in \cup_{r\in R} I_r} \Big( 1- \frac{1}{p^s}\Big)^{-1}
= \sum_{n\in {\Cal S}(R)} \frac{1}{n^s}.
$$

\proclaim{Proposition 4.2}  Let $R$ be any subset of $\{0,\ \ldots, \ J\}$.
Then we have that
$$
\frac 1T \int_{T}^{2T} |\zeta(1+it;R)|^{2z} dt
= \sum_{n\in {\Cal S}(R)} \frac{d_z(n)^2}{n^2}
\Big( 1 + O\Big(\exp\Big( - \frac{\log T}{2(\log_2 T)^4}\Big)  \Big).
$$
\endproclaim

Note that the first part of Theorem 4.1 follows from the case
$R=\{0,1,\ldots, J\}$.  While this is the case of interest for us,
the formulation of Proposition 4.2 is convenient for our proof which is
based on induction on the cardinality of $R$.

\proclaim{Lemma 4.3}  For any prime $p$ we have
$$
\sum_{a=0}^{\infty} \frac{d_z(p^a)^2}{p^{2a}} = I_0\Big(\frac
{2k}{p}\Big) \exp(O(k/p^2)) .
$$
Also
$$
\Big( 1-\frac{\delta}{p}\Big)^{-2k\delta} \ge \sum_{a=0}^{\infty}
\frac{d_z(p^a)^2}{p^{2a}} \geq  \frac{1}{50} \min\Big(1,
\frac{p}{k}\Big) \Big( 1-\frac{\delta}{p}\Big)^{-2k\delta},
$$
so that if $\Cal P$ is any subset of the primes $\leq y$ then, uniformly,
$$
\sum\Sb n\geq 1\\ p|n \implies p\in \Cal P\endSb \frac{ d_z(n)^2}{n^2}
\geq  T^{O(1/\log_2T)} \prod_{p\in \Cal P} \Big(
1-\frac{\delta}{p}\Big)^{-2k\delta}.
$$
\endproclaim
\demo{Proof} Since
$$
\sum_{a=0}^{\infty} \frac{d_z(p^a)^2}{p^{2a}} = \int_0^1 \Big|
1-\frac{e(\theta)}{p}\Big|^{-2z} d\theta = \int_0^1 \exp(O(k/p^2))
\exp\left(2\frac zp \cos (2\pi\theta)\right) d\theta
$$
we obtain the first assertion.  The upper bound in the
second statement follows since $|1-e(\theta)/p|^{-\delta}
\le (1-\delta/p)^{-\delta}$.  When $p>k$ we have
that $(1-\delta/p)^{-2k\delta} \le (1-1/\max(2,k))^{-2k} \le 16$
and so the lower bound follows in this case.
When $p\le k$ consider
only $\theta$ such that $e(\theta)$ lies on the arc $(\delta e^{-ip/(10k)},
\delta e^{ip/(10k)})$.  For such $\theta$ we may check that
$|1-e(\theta)/p|^{-2k\delta} \ge (1-\delta/p)^{-2k\delta} (1-1/(25k))^{k}
\ge \frac 45 (1-\delta/p)^{-2k\delta}$
from which the lower bound in this case follows.

Now
$$
\prod\Sb k<p\leq y \\ p\in \Cal P\endSb \Big( 1-\frac{\delta}{p}\Big)^{-2k\delta} \leq
\exp \Big( O\Big( \sum_{k<p\leq y} \frac {k}p \Big) \Big) 
\ll \left( \frac{\log y}{\log k} \right)^{O(k)} \ll T^{O(1/\log_2T)},
$$
and
$$
\sum\Sb n\geq 1\\ p|n \implies p\in \Cal P\endSb \frac{ d_z(n)^2}{n^2}
> \sum\Sb n\geq 1\\ p|n \implies p\leq k \ \text{and} \ p\in \Cal P\endSb \frac{
d_z(n)^2}{n^2} \geq \prod\Sb p\leq k\\ p\in \Cal P\endSb \frac{p}{50k}  \Big(
1-\frac{\delta}{p}\Big)^{-2k\delta}  ,
$$
which together imply the third assertion by the prime number
theorem.
\enddemo

\proclaim{Lemma 4.4} Suppose $0\le r\le J$ and put $M_0:= T^{\frac 15}$
and $M_r= T^{\frac {1}{5r^2}}$ for $r\ge 1$.
Then we have that
$$
\sum\Sb m \in {\Cal S}(\{r\})\\
m\ge M_r \endSb \frac{2^{\omega(m)}}{m}
\sum_{\ell \in {\Cal S}(\{r\})} \frac{|d_z(m\ell)d_z(\ell)|}{\ell^2}
\le \Big(\sum_{\ell \in {\Cal S}(\{r\})} \frac{d_z(\ell)^2}{\ell^2} \Big)
\exp\Big(-\frac{\log T}{(\log_2 T)^4}\Big).
$$
\endproclaim

\demo{Proof} Denote the left side of the estimate in Lemma 4.4 by
$N_r$ and let 
$$
D_r =\sum_{\ell \in {\Cal S}(\{r\})}
\frac{d_z(\ell)^2}{\ell^2}.
$$ 
For any $1\geq \alpha >0$ we have
$$
\align
N_r
&\le M_r^{-\alpha} \sum\Sb m \in {\Cal S}(\{r\})\endSb
\frac{2^{\omega(m)}}{m^{1-\alpha}}
\sum_{\ell \in {\Cal S}(\{r\})} \frac{|d_z(m\ell)d_z(\ell)|}{\ell^2}  \\
&= M_r^{-\alpha} \prod_{p\in I_r} \Big( \sum_{a=0}^{\infty}
\frac{|d_z(p^a)|^2}{p^{2a}} + 2\sum_{u=1}^{\infty} \frac{1}{p^{u(1-\alpha)}}
\sum_{a=0}^{\infty} \frac{|d_z(p^a)d_z(p^{u+a})|}{p^{2a}}\Big).
\tag{4.1}
\\
\endalign
$$
We record two bounds for the
$p$th term of the product in (4.1):\ Firstly
$$
\align
\sum_{a=0}^{\infty}
\frac{|d_z(p^a)|^2}{p^{2a}} + 2\sum_{u=1}^{\infty} \frac{1}{p^{u(1-\alpha)}}
\sum_{a=0}^{\infty} \frac{|d_z(p^a)d_z(p^{u+a})|}{p^{2a}}
&\le  2\sum_{a=0}^{\infty}
\frac{|d_z(p^a)|}{p^{a(1+\alpha)}} \sum_{u=-a}^{\infty} \frac{|d_z(p^{u+a})|}
{p^{(u+a)(1-\alpha)}}\\
&= 2 \Big(1-\frac{\delta}{p^{1-\alpha}}
\Big)^{-\delta k}\Big(1-\frac{\delta}{p^{1+\alpha}}\Big)^{-\delta k}.
\tag{4.2}\\
\endalign
$$
Secondly, since $|d_z(p^{u+a})|\le |d_z(p^a)||d_z(p^u)|$,
$$
\align
\sum_{ a=0}^{\infty}
\frac{|d_z(p^a)|^2}{p^{2a}} + 2\sum_{u=1}^{\infty} \frac{1}{p^{u(1-\alpha)}}
\sum_{a=0}^{\infty} \frac{|d_z(p^a)d_z(p^{u+a})|}{p^{2a}}
&\le \sum_{a=0}^{\infty} \frac{|d_z(p^a)|^2}{p^{2a}} \Big( 1
+2\sum_{u=1}^{\infty} \frac{|d_z(p^u)|}{p^{u(1-\alpha)}}\Big) \\
&\le \sum_{a=0}^{\infty} \frac{|d_z(p^a)|^2}{p^{2a}} \Big(
2\Big(1-\frac{\delta}{p^{1-\alpha}}\Big)^{-\delta k} -1\Big). \tag{4.3}
\\
\endalign
$$

Now consider the case $r=0$ and note that $k\le p$ for all $p\in I_0$.
Here we use the bound (4.3) in (4.1).  We choose $\alpha= 1/(10\log_2 T)$ and
note that for $p\in I_0$, $2(1-\delta/p^{1-\alpha})^{-\delta k} -1
\le 2(1-e^{1/9}/p)^{-k}-1 \le e^{4k/p}$.  Hence we get that
$$
\align
N_0 &\le D_0 \exp\Big( -\frac{\log M_0}{10 \log_2 T} + 4k\sum_{k\le p\le  y}
\frac 1p \Big)
\le D_0 \exp\Big( -\frac{\log M_0}{10 \log_2 T} + \frac{4k}{\log k}
\sum_{k\le p\le y}\frac {\log p}{p}\Big). \\
\endalign
$$
Now $\sum_{k\le p\le y} \log p/p\le \log(25y/k)$ (see Theorem I.1.7
of Tenenbaum [7]) and recall that $k\le \log T/(e^{10} \log (y/\log
T))$ and that $M_0 =T^{1/5}$.  The bound in the lemma then follows
in this case.

Suppose now that $r\ge 1$ so that $p\le k$ for all $p\in I_r$.  Here
we use the bound (4.2) in (4.1).   We take $\alpha
= 1/(10\cdot 2^{r/2} \log (ek))$ and note that for $p\le k$,
$$
\align
\Big( 1-\frac{\delta}{p^{1-\alpha}}\Big)^{-\delta}
\Big(1-\frac{\delta}{p^{1+\alpha}}\Big)^{-\delta}
\Big(1-\frac{\delta}{p}\Big)^{2\delta}
&\le \Big( 1-\frac{p(p^\alpha +p^{-\alpha} -2)}{(p-1)^2}\Big)^{-1} \\
&\le \exp\Big( \frac{\log^2 p}{10\cdot 2^r p \log^2 (ek)}\Big).
\\
\endalign
$$
Using also the lower bound in Lemma 4.3 we obtain that
$$
N_r \le D_r \exp\Big( -\frac{\log M_r}{10 \cdot 2^{r/2} \log (ek)} +
\sum_{p\in I_r} \Big( \log \frac{100 k}{p} + \frac{k\log p}{10 \cdot
2^r p \log (ek)}\Big)\Big). \tag{4.4}
$$
If $1\le r\le J-1$ then we deduce that
$$
\align
N_r &\le D_r \exp\Big(-\frac{\log M_r}{10 \cdot 2^{r/2} \log (ek)} +
\sum_{k/2^r \le p \le k/2^{r-1}} (r+5)\Big) \\
&\le D_r \exp\Big(-\frac{\log M_r}{10 \cdot 2^{r/2} \log (ek)} +
\frac{8(r+5)k}{2^r \log (ek)}\Big)
\\
\endalign
$$
and since $\log M_r = (\log T)/(5r^2)$ this gives
$N_r\le D_r \exp(-\log T/(\log_2 T)^4)$ for large $T$.
If $r=J$ and $k\le 10^6$ then the Lemma follows at once
from (4.4).  If $r=J$ and $k>10^6$ then (4.4) gives that
$$
\align
N_r&\le D_r \exp\Big( -\frac{\log M_J}{10\cdot 2^{J/2} \log (ek)}
+ \sum_{p\le k/(\log k)^4} \Big(\log \frac{100k}{p}
+ \frac{k\log p}{10 \cdot 2^J p \log(ek)}\Big)\Big)
\\
&\le D_r \exp\Big(-\frac{\log M_J}{10\cdot 2^{J/2} \log (ek)} +O\Big(\frac{
\log T }{(\log_2 T)^4}\Big) \Big),\\
\endalign
$$
which proves the Lemma in this case.

\enddemo

\demo{Proof of Proposition 4.2} We prove Proposition 4.2 by induction
on the cardinality of $R$.  The case when $R =\emptyset$ is clear
and suppose the Proposition holds for all proper subsets of $R$.
We expand
$$
|\zeta(1+it;R)|^{2z} = \sum\Sb m_r, n_r \in {\Cal S}(\{r\})  \\
\text{for all} \ r\in R\endSb \ \prod_{r\in R} \left(
\frac{d_z(m_r) d_z( n_r)}{ m_r n_r } \right) \Big(\frac{\prod_{r\in
R} m_r }{\prod_{r\in R} n_r}\Big)^{it}.
$$
Set $u_r = m_r n_r/(m_r,n_r)^2$.
Using inclusion-exclusion we decompose the above as
$$
\sum\Sb m_r, n_r \in {\Cal S}(\{r\}),\ \text{and}  \\
u_r\le M_r \ \text{for all} \ r\in R
\endSb + \sum\Sb W\subset R\\ W\neq \emptyset \endSb
(-1)^{|W|-1} \sum\Sb m_r, n_r \in {\Cal S}(\{r\})   \\
\text{for all} \ r\in R, \ \text{and} \\ u_w > M_w \ \text{for all} \ w\in W
\endSb  \tag{4.5}
$$
with $M_w$ as in Lemma 4.4.

First let us consider the contribution of the first sum in (4.5).  This
gives
$$
\sum\Sb m_r, n_r \in {\Cal S}(\{r\}),\ \text{and}  \\
u_r\le M_r \ \text{for all} \ r\in R\endSb
\ \prod_{r\in R} \left( \frac{d_z(m_r) d_z( n_r)}{ m_r n_r } \right)
\frac 1T \int_{T}^{2T} \Big( \prod_{r\in R} \frac{m_r }{n_r
}\Big)^{it} dt. \tag{4.6}
$$
If we reduce $\prod_{r\in R} m_r/n_r$ to lowest terms then both
the numerator and denominator would be bounded by $\prod_{r} u_r
\le \prod_{r\in R} M_r \le T^{\frac{(1+\pi^2/6)}{5}}\le T^{\frac 35}$.
Thus if $\prod_{r\in R} m_r/n_r \neq 1$ then
$$
\frac{1}{T}\int_{T}^{2T}
\Big( \frac{\prod_{r\in R} m_r }{\prod_{r\in R} n_r }\Big)^{it} dt
\ll \frac{1}{T |\log \prod_{r} m_r/n_r|} \ll T^{-\frac 25}.
$$
Hence we obtain that the expression in (4.6) equals
$$
\sum\Sb m_r=n_r \in {\Cal S}(\{r\}) \\ \text{for all} \ r\in R\endSb \ \prod_{r\in R}
\left( \frac{d_z(m_r)}{m_r} \right)^2  +
O\Big( T^{-\frac 25} \sum\Sb  m_r, n_r \in {\Cal S}(\{r\}) \\ \text{for all} \ r\in R \endSb 
\ \prod_{r\in R} \left( \frac{|d_z(m_r) d_z( n_r)|}{ m_r n_r } \right)\Big).
$$
The main term above is $\sum_{n\in {\Cal S}(R)} d_z(n)^2/n^2$.
The error term is $\ll T^{-\frac 25} \prod_{p\in \cup_{r\in R} I_r }
(1-\delta/p)^{-2k\delta}$ and using the lower bound of
Lemma 4.3 this is $\ll T^{-\frac 13} \sum_{n\in {\Cal S}(R)} d_z(n)^2/n^2$.
Thus the contribution of the first term in (4.5) is
$$
(1+O(T^{-\frac 13})) \sum_{n\in {\Cal S}(R)} \frac{d_z(n)^2}{n^2}.
\tag{4.7}
$$

Now we consider the contribution of the second term in (4.5).
This gives
$$
\align \sum\Sb W \subset R \\ W\neq\emptyset \endSb (-1)^{|W|-1}
&\sum\Sb m_w, n_w \in {\Cal S}(\{w\}),\ \text{and}  \\
u_w>M_w \ \text{for all} \ w\in W \endSb \ \prod_{w\in W}
\left( \frac{d_z(m_w) d_z( n_w)}{ m_w n_w } \right)\\
&\times \frac{1}{T} \int_T^{2T} \Big(\frac{\prod_{w\in W}
m_w}{\prod_{w\in W} n_w}
\Big)^{it} |\zeta(1+it;R-W)|^{2z} dt,\\
\endalign
$$
which is bounded in magnitude by
$$
\sum\Sb W \subset R \\ W\neq\emptyset \endSb \ \sum\Sb m_w, n_w \in {\Cal S}(\{w\}),\ \text{and}  \\
u_w>M_w \ \text{for all} \ w\in W \endSb \ \prod_{w\in W} \left( \frac{|d_z(m_w) d_z(
n_w)|}{ m_w n_w } \right) \frac{1}{T} \int_T^{2T}
|\zeta(1+it;R-W)|^{2z} dt.
$$
By the induction hypothesis we see that
$$
\frac 1T \int_{T}^{2T} |\zeta(1+it;R-W)|^{2z} dt \ll \sum_{n \in
{\Cal S}(R-W)} \frac{d_z(n)^2}{n^2},
$$
while from Lemma 4.4 (with $m=u_w$ and $\ell=(m_w,n_w)$ so that $
d_z(m\ell)d_z(\ell)=d_z(m_w)d_z(n_w)$; and note that the number of pairs 
$m_w,n_w$ which give rise to a given pair $\ell, m$ is exactly $2^{\omega(m)}$) we deduce that
$$
\sum\Sb m_w, n_w \in {\Cal S}(\{w\}) \\ u_w >M_w \endSb
\frac{|d_z(m_w)d_z(n_w)|}{m_w n_w} \le \sum_{n\in {\Cal S}(\{w\})}
\frac{d_z(n)^2}{n^2} \exp\Big( -\frac{\log T}{(\log_2 T)^4}\Big).
$$
From these estimates it follows that the contribution of the
second term in (4.5) is
$$
\ll  |R| \sum_{n\in {\Cal S}(R)} \frac{d_z(n)^2}{n^2}
\exp\Big( -\frac{\log T}{(\log_2 T)^4}\Big).
$$
Combining this with (4.7) we obtain Proposition 4.2.

\enddemo

\demo{Proof of Theorem 4.1} In view of Proposition 4.2 it remains
only to prove that
$$
\sum\Sb n=1\\ p|n\implies p\le y\endSb^{\infty} \frac{d_z(n)^2}{n^2}
= \prod_{p\le k} \Big(1-\frac{\delta}{p} \Big)^{-2k\delta}
\exp\Big(\frac{2k}{\log k} \Big(C +O\Big(\frac{k}{y} +\frac{1}{\log
k}\Big)\Big)\Big). \tag{4.8}
$$
Using the first part of Lemma 4.3 for $p\ge \sqrt{k}$ and
the second part for $p<\sqrt{k}$ we see that
$$
\sum\Sb n=1\\ p|n\implies p\le y\endSb^{\infty}
\frac{d_z(n)^2}{n^2}  = \prod_{p<\sqrt{k}} \Big(1-\frac \delta{p}\Big)^{-2k
\delta} \prod_{\sqrt {k} \le p \le y} I_0\Big(\frac {2k}{p}\Big)
\exp(O(\sqrt{k})).
$$
Since $\log I_0(t)= O(t^2)$ for $0\le t\le 2$ we
have by the prime number theorem and partial summation that
$$
\align
\sum_{k\le p\le y} \log I_0\Big(\frac{2k}{p} \Big)
&=\frac{2k}{\log k} \int_{2k/y}^{2} \log I_0(t) \frac{dt}{t^2}
+ O\Big(\frac{k}{\log^2 k}\Big)
\\
&= \frac{2k}{\log k} \int_0^2 \log I_0(t) \frac{dt}{t^2} +
O\Big( \frac{k^2}{y\log k} + \frac{k}{\log^2 k}\Big).
\\
\endalign
$$
Since $\log I_0(t) = t+ O(\log t)$ for $t\ge 2$ we obtain by the
prime number theorem and partial summation that
$$
\sum_{\sqrt{k}\le p\le k} \Big(\log I_0\Big(\frac{2k}{p}\Big) +
2k\delta \log \Big(1-\frac \delta{p}\Big)\Big) = \frac{2k}{\log k}
\int_2^{\infty} (\log I_0(t) -t)\frac{dt}{t^2}
+O\Big(\frac{k}{\log^2 k}\Big).
$$
These estimates prove (4.8) and so Theorem 4.1 follows.

\enddemo

\head 5.  Proof of Theorem 1 \endhead

\noindent Let $\log T (\log_2 T)^4 \ge y\ge e^2 \log T$, and let
 $T\Phi_T(\tau;y)$ denote the measure of points
$t\in [T,2T]$ for which $|\zeta(1+it;y)| \ge e^{\gamma} \tau$.
Taking $z=k$ for an integer $3\le k \le \log T/(e^{10} \log (y/\log T))$
in Theorem 4.1 and using Mertens' theorem
$\prod_{p\le k} (1-1/p)^{-1} = e^{\gamma} \log k
+ O(1/\log^2 k)$ we get that
$$
\align
2k\int_0^\infty \Phi_T(t;y) t^{2k-1} dt
&= \frac{1}{T} \int_{T}^{2T} e^{-2k\gamma} |\zeta(1+it;y)|^{2k} dt \\
&= ( \log k)^{2k}
\exp\left(\frac{2k}{\log k}\Big(C +O\Big(\frac{k}{y} +
\frac{1}{\log k}\Big)\Big)
\right).
\tag{5.1}\\
\endalign
$$
Now $\int_{0}^{\infty} \Phi_T(t;y) dt = e^{-\gamma} (1/T)
\int_{T}^{2T}|\zeta(1+it;y)| dt \leq e^{-\gamma} ((1/T)
\int_{T}^{2T}|\zeta(1+it;y)|^4 dt)^{1/4}\ll 1$ by Theorem 4.1; so,
by H{\" o}lder's inequality,
$$
\int_{0}^{\infty} \Phi_T(t;y) t^a dt \leq \Big(\int_{0}^{\infty}
\Phi_T(t;y)  dt\Big)^{1-a/b} \Big(\int_{0}^{\infty} \Phi_T(t;y)  t^b
dt\Big)^{a/b} \ll \Big(\int_{0}^{\infty} \Phi_T(t;y)  t^b
dt\Big)^{a/b}
$$
for $a<b$.
 While (5.1) at present holds only for integer values of
$k$, we may interpolate to non-integer value $\kappa\in (k-1,k)$ by
taking $a=2k-3 ,\ b=2\kappa-1$ and then $a=2\kappa-1 ,\ b=2k-1$ in
the last inequality to obtain
$$
\Big( \int_{0}^{\infty} \Phi_T(t;y) t^{2k -3}dt \Big)^{\frac{2\kappa
-1}{2k-3}} \ll \int_0^{\infty} \Phi_T(t;y)t^{2\kappa -1} dt \ll
\Big(\int_0^{\infty} \Phi_T(t;y)t^{2k-1} dt
\Big)^{\frac{2\kappa-1}{2k-1}},
$$
and so we get (5.1) for $\kappa$ by substituting (5.1) for $k-1$ and
$k$ into this equation.


Suppose $1\ll \tau \le \log_2 T -20 - \log_2 (y/\log T)$ and select
$\kappa=\kappa_\tau$ such that $\log \kappa = \tau - 1 - C$.
Let $\epsilon>0$ be a bounded parameter to be fixed shortly and put
$K=\kappa e^{\epsilon}$.
Observe that
$$
\align
2k\int_{\tau+\epsilon}^\infty \Phi_T(t;y) t^{2\kappa-1} dt
&\leq 2k {(\tau+\epsilon)}^{2\kappa-2K} \int_{\tau+\epsilon}^\infty
\Phi_T(t;y) t^{2K-1} dt
\\
& \leq (\tau+\epsilon)^{2\kappa(1-e^{\epsilon})}
\Big(2K \int_0^\infty \Phi_T(t;y) t^{2K-1} dt\Big).
\\
\endalign
$$
Using (5.1) we deduce that
$$
\align
&2K\int_0^{\infty} \Phi_T(t;y) t^{2K-1} dt
\\
&= \Big( (\log \kappa +\epsilon)
\exp\Big(\frac{C}{\log \kappa} \Big(1+O\Big(\frac{1}{\log \kappa}+
\frac{\kappa}{y} \Big) \Big)
\Big)\Big)^{2K}\\
&= \exp\Big( \frac{2\kappa(\epsilon e^{\epsilon} + C(e^{\epsilon}-1))}
{\log \kappa}+O\Big(\frac{\kappa}
{\log^2 \kappa}+\frac{\kappa^2}{y\log \kappa}\Big)\Big)
(\log \kappa)^{2\kappa(e^{\epsilon}-1)}
\int_0^{\infty} \Phi_T(t;y) t^{2\kappa-1} dt.
\\
\endalign
$$
We conclude from the above that
$$
2\kappa
\int_{\tau+\epsilon}^{\infty} \Phi_T(t;y) t^{2\kappa -1} dt
= \exp\Big(\frac{2\kappa}
{\log \kappa} (1+\epsilon-e^{\epsilon}) + O\Big(\frac{\kappa}{\log^2 \kappa}
+\frac{\kappa^2}{y\log \kappa}\Big)\Big)
\int_0^{\infty} \Phi_T(t;y) t^{2\kappa -1} dt.
$$
Choose  $\epsilon = c(1/\tau + (\log T)/y)^{\frac 12}$ for a suitable 
constant $c>0$, so that for large $\tau$ (and hence large $\kappa$),
$$
\int_{\tau+ \epsilon}^{\infty} \Phi_T(t;y) t^{2\kappa -1} dt
\le \frac{1}{100} \int_0^{\infty} \Phi_T(t;y) t^{2\kappa -1} dt,
$$
say.
A similar argument reveals that
$$
\int_0^{\tau-\epsilon} \Phi_T(t;y) t^{2\kappa -1} dt
\le \frac{1}{100} \int_0^{\infty} \Phi_T(t;y) t^{2\kappa -1} dt.
$$

Combining these two assertions with (5.1) for $\kappa$ we obtain
$$
\int_{\tau-\epsilon}^{\tau+\epsilon}
\Phi_T(t;y) t^{2\kappa-1} dt = (\log \kappa)^{2\kappa}
\exp\Big(
\frac{2\kappa C}{\log \kappa} (1+O(\epsilon^2)) \Big).
$$
Since $\Phi_T$ is a non-increasing function we deduce that the left
side above is
$$
\ge \Phi_T(\tau+\epsilon;y) \tau^{2\kappa} \exp(O(\kappa \epsilon/\tau)),
\qquad
\text{and}
\qquad
\le \Phi_T(\tau-\epsilon;y) \tau^{2\kappa} \exp(O(\kappa\epsilon/\tau))
.
$$
It follows that
$$
\Phi_T(\tau+\epsilon;y) \le \exp\Big( -(2
+ O(\epsilon)) \frac{e^{ \tau-1-C}}{\tau}
\Big)\le \Phi_T(\tau-\epsilon;y),
$$
and hence that uniformly in $\tau \le \log_2 T -20 - \log_2 (y/\log T)$
we have
$$
\Phi_T(\tau;y) = \exp\Big( - \frac{2e^{\tau -1-C}}{\tau} (1+O(\epsilon)) \Big)\Big). \tag{5.2}
$$
From Proposition 3.2 we know that $\Phi_T(\tau)
= \Phi_T(\tau+O(\epsilon);y)
+ O(\exp(-\log T/50\log_2 T))$ for $\tau\ll \log_2T$;
and so from (5.2) we deduce
that uniformly in $\tau \le \log_2 T -20 -\log (y/\log T)$ we
have 
$$
\Phi_T(\tau) =\exp\Big( - \frac{2e^{\tau -1-C}}{\tau} (1+O(\epsilon))\Big) +O\Big(\exp\Big(-\frac{\log T}{50\log_2 T}\Big)
\Big).
$$
Taking $y=\min( \tau \log T, (\log^2 T)/e^{10+\tau})$ above we
easily obtain Theorem 1 for $\Phi_T$.  The argument for $\Psi_T$ is
analogous, using $z=-k$ in Theorem 4.1.

\bigskip

One finds, using the first part of Lemma 4.3  and
the observation that $\log I_0(2k/p) \ll k^2/p^2$ for $p>k$, that
$$ \align
\Bbb E( |L(1,X)|^{2z} )  &= \sum_{n\geq 1} \frac{d_z(n)^2}{n^2} =\sum\Sb n=1\\ p|n \implies p\le y\endSb^{\infty}
\frac{d_z(n)^2}{n^2} \exp\Big(  O\Big(\frac{k^2}{y\log y}  \Big)\Big)
\\
&=\prod_{p\le k} \Big(1-\frac{\delta}{p}\Big)^{-2k\delta}
\exp\Big(\frac{2k}{\log k}\Big(C+ O\Big(\frac{k}{y} +\frac{1}{\log k}\Big)
\Big)\Big),\\
\endalign
$$
the last line following as in the proof of Theorem 4.1. With this estimate we
can proceed precisely as in the proof of Theorem 1 to obtain the
analagous estimate.

 \head 6.  Large values of $|\zeta(1+it)|$: Proof of Theorem 2\endhead

\noindent Let $T$ be large and put $y= \log T \log_2 T/(4B \log_3
T)$ for some $B\geq 5$, and $\delta = 1/[\log_2 T]^4$.  Let
$\parallel z\parallel$ denote the distance of $z$ from the nearest
integer.

\proclaim{Lemma 6.1}  For any real $t_0$  there is a positive
integer $m\le T^{\frac 1B}$ such that for 
each prime $p\le y$ we have $\parallel (mt_0 \log p)/2\pi
\parallel \le \delta$.
\endproclaim

\demo{Proof}  This follows from Dirichlet's theorem on Diophantine
approximation (see for example \S 8.2 of  [8]) since $1/\delta$ is
an integer and  $(1/\delta)^{\pi (y)} \le T^{\frac 1B}$, by the prime number theorem.
\enddemo

\proclaim{Lemma 6.2}  For any real $t_1$ there is a positive integer
$n \le [\log_2 T]^2$ for which
$$
\text{\rm Re }\sum_{y\le p\le \exp((\log T)^{10})} \frac{1}{p^{1+int_1}}
\ge - \frac{10}{\log_2 T}.
$$
\endproclaim
\demo{Proof} Let $K(x)=\max(0,1-|x|)$ and note that
$\sum_{l=-L}^L K(l/L)e^{ilt}$ (the Fejer kernel) is non-negative for
all positive integers $L$ and all $t$.  It follows
therefore that
$$
\sum_{j=- [\log_2 T]^2}^{ [\log_2 T]^2} K\Big( \frac{j}{[\log_2 T]^2}
\Big) \sum_{y\le p\le \exp((\log T)^{10})} \frac{1}{p^{1+ijt_1}}
\ge 0.
$$
Hence we obtain that
$$
\align
\text{Re } \sum_{j=1}^{ [\log_2 T]^2} K\Big( \frac{j}{[\log_2 T]^2}
\Big) \sum_{y\le p\le \exp((\log T)^{10})} \frac{1}{p^{1+ijt_1}}
&\ge - \frac{1}{2} \sum_{y\le p\le \exp((\log T)^{10})} \frac{1}{p} \\
&\ge - 5\log_2 T.\\
\endalign
$$
The Lemma follows at once.
\enddemo

\demo{Proof of Theorem 2}
 For $T^{\frac 1{10}} \le |t| \le T$ one has
$$
\log \zeta(1+it)
=- \sum_{p \le \exp((\log T)^{10})} \log \Big(1-\frac{1}{p^{1+it}}
\Big) + O\Big(\frac{1}{\log T}\Big).
$$
(One can prove this, arguing as in the proof of the prime number
theorem, by noting that $(1/2i\pi) \int_{(c)} \log \zeta(1+it+w) (x^w/w)dw$
with $x=\exp((\log T)^{10})$ and $c>0$ gives the main term of the
right side by Perron's formula, and by shifting the contour to the
left of $0$, but enclosing a region free of zeros of $\zeta(s)$, we
get residue $\log \zeta(1+it)$ from the simple pole at $w=0$, and
the error term from the remaining integral.)

Combining Lemmas 6.1 and 6.2 (with $t_1=mt_0$) we see that for any
$t_0\in [T^{1/10},T]$ there exists an integer $\ell$ (where $\ell=mn$) 
with $1\le \ell\le T^{\frac 1B} [\log_2 T]^2$ such that $\parallel (\ell t_0 \log
p)/2\pi\parallel \le 1/[\log_2 T]^2$ for each prime $p\leq y$, and
such that
$$
\text{Re }\sum_{y\le p\le \exp((\log T)^{10})} \frac{1}{p^{1+i\ell t_0}}
\ge -\frac{10}{\log_2 T}.
$$
We deduce therefore that
$$
\align
|\zeta(1+i\ell t_0)| &\ge
\prod_{p\le y} \Big (1-\frac{1}{p} + O\Big(\frac{1}{p(\log_2 T)^2}\Big)
\Big)^{-1} \Big(1+O\Big(\frac{1}{\log_2 T}\Big) \Big) \\
&\ge e^{\gamma} (\log_2 T +\log_3 T -\log_4 T - \log A+O(1)),\\
\endalign
$$
using the prime number theorem, where $A=1/(2/B+3\log_2T/\log T)$.

We use the above procedure with $t_0 = T_0, \ T_0 +1, \ T_0 +2,
\ldots, \ T_0 + U_0$ where $T_0=[T^{1-1/B}/3[\log_2T]^2]$ and
$U_0=[T^{1-2/B}/7[\log_2T]^4]$ . Let $\ell_i$ be as above so
$\ell_i\leq T^{1/B}[\log_2T]^2$ and thus $\tau_i=\ell_i (T_0+i) \leq
T/2$. We claim that $|\tau_i-\tau_j|\geq 1$ if $i\ne j$ for if not
then evidently $\ell_i\ne \ell_j$ (else $1\leq |(T_0+j)-(T_0+i)|=|\tau_j-\tau_i|/\ell_i<1$), 
so that
$$
T_0\leq |(\ell_i-\ell_j)T_0| \leq |\tau_i-\tau_j|
+|i\ell_i-j\ell_j|< 1+U_0T^{1/B}[\log_2T]^2,
$$
which is false. Now each $|\zeta(1+i\tau_j)| \ge e^{\gamma} (\log_2
T +\log_3 T -\log_4 T-\log A+O(1))$. Since $|\zeta^{\prime}(1+it)|
\ll \log^2 T$ for $1\le |t|\le T$ we see that for any $|\alpha|\le
1/\log^2 T$ we have that $|\zeta(1+i\tau_j + i\alpha)| =
|\zeta(1+i\tau_j)| +O(\alpha \log^2 T) = |\zeta(1+i\tau_j)| +O(1)$.
Thus the measure of $t \in [0,T]$ with $|\zeta(1+it)| \ge e^{\gamma}
(\log_2 T +\log_3 T -\log_4 T-\log A+O(1))$ is at least $2U_0/\log^2
T$, proving Theorem 2.
\enddemo

\head 7. The analogous results for $L$-functions at $1$ \endhead

\noindent By analogous methods one can prove:

\proclaim{Theorem 3}  Let $q$ be a large prime. \hfill \break (i)
The proportion of characters $\chi \pmod q$ for which $|L(1,\chi)|
>e^{\gamma}\tau$ is
$$
\exp\Big(-\frac{2e^{\tau-C-1}}{\tau} \Big(1 +
O\Big(\frac{1}{\tau^{\frac 12}} +\Big(\frac{e^{\tau}}{\log
q}\Big)^{\frac 12}
 \Big)\Big)\Big) ,  \tag{7.1}
$$
uniformly in the range $1\ll \tau \le \log_2 q -20$. The same
asymptotic also holds for the proportion of characters $\chi \pmod
q$ for which $|L(1,\chi)| <{\pi^2}/{6e^{\gamma}\tau}$. \hfill \break
(ii) There are at least $q^{1-1/A}$ characters $\chi \pmod q$  such
that
$$
|L(1,\chi)| \ge e^{\gamma}(\log_2 q + \log_3 q- \log_4 q- \log A
+O(1)),
$$
 for any $A\geq 10$.
\endproclaim







If, in addition, we vary over all characters $\chi \pmod q$ 
and all primes $Q\le q\le 2Q$, then we can get a good estimate for the
distribution function of $|L(1,\chi)|$ in almost the
entire viable range.
Thus we may prove that the proportion of $|L(1,\chi)|\ge
e^{\gamma} \tau$ is (7.1) for the range $1\le \tau \le \log_2 Q +
\log_3 Q -100$, but now with the error term ``$(e^\tau/(\log Q
\log_2Q))^{\frac 12}$'' in place of ``$(e^\tau/\log q)^{\frac 12}$''
(and a corresponding result holds for $1/|(6/\pi^2)L(1,\chi)|$).

The broad outline of the proof is the same, though now replacing
$\log T$ by $\log Q \log_2 Q$, so that $ \log Q(\log_2 Q)^4 \ge y
\ge e^2 \log Q \log_2 Q$ and the range for $k$ becomes 
$2\le k \le \log Q \log_2 Q/(e^{10}\log
(y/(\log Q \log_2 Q)))$. The result follows easily from the
following analogy to Theorem 4.1,
$$
\frac 1{\pi(Q)} \sum_{q\le Q} \frac{1}{\phi(q) } \sum_{\chi\pmod q}
|L(1,\chi;y)|^{2z} =\prod_{p\le k} \Big( 1-\frac{\delta}{p}\Big)^{-2
k\delta} \exp\Big( \frac{2k}{\log k} \Big(C_1 +O\Big( \frac{k}{y}
+\frac{1}{\log k}\Big)\Big)\Big) ,
$$
and an appropriate development of Lemma 4.3, where
$L(1,\chi;y):=\prod_{p\le y} (1-\chi(p)/p)^{-1}$. The above estimate,
though, is proved rather more easily than Theorem 4.1. Since
$L(1,\chi;y)^z=\sum_{n\in {\Cal S}(y)} d_z(n)\chi(n)/n$, and
$L(1,\cbar;y)^{z} =\sum_{m\in {\Cal S}(y)} d_z(m)\cbar(m)/m$ where
$S(y)$ is the set of integers all of whose prime factors are $\leq
y$ , the left side of this equation equals
$$
\sum\Sb m,n \in {\Cal S}(y)\endSb \frac{d_z(m) d_z(n)}{mn} \Big\{
\frac{1}{\pi(Q)} \sum_{q\le Q} \frac{1}{\phi(q) } \sum_{\chi \pmod
q} \chi(m)\cbar(n) \Big\}.
$$
The term in $\{\}$ equals $1-\#\{ q\leq Q:\ q|mn\}/\pi(Q)$ if $m=n$, and
is $\leq \#\{ q\leq Q:\ q|m-n\}/\pi(Q)$ if $m\ne n$. Therefore our sum is
$$
\sum_{n\in {\Cal S}(y)} \frac{d_{z}(n)^2}{n^2} + O\Big(
\frac{1}{\pi(Q)} \Big( \sum_{m\in {\Cal S(y)}} \frac{|d_z(m)| \log
2m}{m} \Big)^2 \Big) .
$$
Now $\log 2n \ll k^2 + n^{1/k}$ so that
$$
\align \sum_{n\in {\Cal S}(y)} \frac{|d_z(n)|}{n} \log 2n  &\ll k^2
\prod_{p\le y} \Big(1-\frac{\delta}{p}\Big)^{-\delta k} +
\prod_{p\le y} \Big( 1-\frac{\delta}{p^{1-1/k}} \Big)^{-\delta k}
\\
&\ll \prod_{p\le y} \Big(1-\frac{\delta}{p}\Big)^{-\delta k} \Big(
k^2 + \exp\Big( O\Big( k \sum_{p\le y} \frac{p^{1/k} -1}{p}
\Big)\Big) \Big)
\\
&\ll
 (\log Q)^{O(1)} \prod_{p\le y}
\Big(1-\frac{\delta}{p}\Big)^{-\delta k} ,
\\
\endalign
$$
and the claimed estimate follows from Lemma 4.3.

\enddocument